\documentclass{amsart}

\usepackage{amsmath,amsthm,amssymb,amsfonts,enumerate}

\numberwithin{equation}{section} 
\theoremstyle{plain}
\newtheorem{theo+}           {Theorem}      [section]
\newtheorem{prop+}  [theo+]  {Proposition}
\newtheorem{coro+}  [theo+]  {Corollary}
\newtheorem{lemm+}  [theo+]  {Lemma}
\newtheorem{defi+}  [theo+]  {Definition}

\theoremstyle{definition}
\newtheorem{exam+}  [theo+]  {Example}
\newtheorem{rema+}  [theo+]  {Remark}

\newenvironment{theorem}{\begin{theo+}}{\end{theo+}}
\newenvironment{proposition}{\begin{prop+}}{\end{prop+}}
\newenvironment{corollary}{\begin{coro+}}{\end{coro+}}
\newenvironment{lemma}{\begin{lemm+}}{\end{lemm+}}

\newenvironment{remark}{\begin{rema+}}{\end{rema+}}

\input diagrams   
\newarrow{Corresponds}<--->

\newcommand{\al}{\alpha}
\newcommand{\be}{\beta}       
\newcommand{\ga}{\gamma}  
\newcommand{\de}{\delta}
\newcommand{\qb}[2]{\genfrac{[}{]}{0pt}{}{#1}{#2}_{q}}
\newcommand{\id}{\operatorname{id}}
\begin{document}

\baselineskip 18pt
\larger[2]
\title
[An elementary approach to  $6j$-symbols]
{An elementary approach to $6j$-symbols\\ (classical, quantum,
rational,\\  trigonometric, and elliptic)} 
\author{Hjalmar Rosengren}
\address
{Department of Mathematics\\ Chalmers University of Technology and G\"oteborg
 University\\SE-412~96 G\"oteborg, Sweden}
\email{hjalmar@math.chalmers.se}
\urladdr{http://www.math.chalmers.se/{\textasciitilde}hjalmar}

\dedicatory{\large Dedicated to Richard Askey}

\begin{abstract}
Elliptic $6j$-symbols first appeared in connection with solvable
models of statistical mechanics. They include many interesting limit
cases, such as quantum $6j$-symbols (or $q$-Racah polynomials) and
Wilson's biorthogonal ${}_{10}W_9$ functions. We give an elementary
construction of elliptic $6j$-symbols, which immediately implies several
of their main properties. As a consequence, we obtain a new algebraic 
interpretation of elliptic $6j$-symbols in terms of
  Sklyanin algebra representations.
\end{abstract}

\maketitle        

\section{Introduction}  

The classical $6j$-symbols were introduced by Racah and Wigner in the
early  1940's \cite{ra,w}. Though they appeared in the  context of
quantum mechanics,  they are natural objects in the representation theory of
$\mathrm{SL}(2)$ that
can be introduced from purely mathematical
considerations. Wilson \cite{wi}  realized that $6j$-symbols are
orthogonal polynomials, and that they
  generalize many classical systems such as
Krawtchouk and Jacobi polynomials. This led Askey and Wilson to
introduce the more general $q$-Racah polynomials \cite{aw}.

The $q$-Racah polynomials belong to the class of basic (or $q$-)
hypergeometric series \cite{gr}. Since
the 1980's,  there has been a considerable increase of
interest in this classical subject. One  reason for this is 
relations to solvable models in statistical mechanics,
 and to the related algebraic structures known as quantum groups. 

Kirillov and Reshetikhin \cite{kr} found that  $q$-Racah
polynomials appear as $6j$-symbols of the
$\mathrm{SL}(2)$ quantum group, or 
\emph{quantum $6j$-symbols}. We mention that in the introduction to the 
standard reference \cite{cp}, three major applications of quantum groups to
other fields of mathematics are highlighted. 
 For at least two of
these, namely,
invariants of links and three-manifolds \cite{t}, and the relation to
affine Lie algebras and conformal field theory \cite{e}, quantum
$6j$-symbols play a decisive role. 

The $q$-Racah polynomials form, together with the  closely related 
 Askey--Wilson polynomials,
 the top level of the Askey Scheme of ($q$-)hypergeometric orthogonal
 polynomials \cite{ks}. One reason for viewing this scheme as
 complete is  Leonard's theorem \cite{l}, saying that
any finite system of orthogonal polynomials with polynomial duals is
a special or degenerate case of the $q$-Racah polynomials. However, if
 one is willing to pass from orthogonal polynomials to biorthogonal
 rational functions,  natural  extensions of the Askey
 Scheme do exist.

One such extension was found by Wilson
  \cite{w2}, who constructed a system of
 biorthogonal rational functions given by 
  ${}_{10}\phi_9$ (or, more precisely, ${}_{10}W_9$) 
basic hypergeometric series.  
These form a generalization of 
$q$-Racah polynomials that seems very natural from the viewpoint of
special functions; see also \cite{rs}.

Another indication that natural generalizations of quantum
$6j$-symbols exist
 came from statistical mechanics. The solvable models
that lead to standard quantum groups and quantum $6j$-symbols appear
there as degenerate cases. Typically, the most
general case of the models  involve
elliptic functions. In the  1980's, Date et al.\ \cite{d2,d3} 
 developped a
fusion procedure for constructing generalized $6j$-symbols from
 $R$-matrices of face models. When applied to 
Baxter's eight-vertex SOS model \cite{abf,b2}, 
this leads to \emph{elliptic $6j$-symbols}, which include
quantum $6j$-symbols as a degenerate case.
  However, no identification of these
objects with biorthogonal rational 
functions was obtained, nor was their nature as generalized
hypergeometric sums emphasized.

In the latter direction,  Frenkel and Turaev \cite{ft1} found that
the \emph{trigonometric} limit case of  elliptic $6j$-symbols  can
be written as ${}_{10}W_9$-series. A further limit transition
 gives \emph{rational} $6j$-symbols. 
Moreover, in \cite{ft2} it was found that  general
elliptic $6j$-symbols may be expressed as elliptic, or modular,
hypergeometric series, a completely new class of special functions. In
spite of their intriguing properties, including  close relations to
elliptic functions and modular forms, such series were never
considered in ``classical'' mathematics, but needed physics for their
discovery.  We refer to \cite{gr2} for an
introduction to the subject, with further references.

Frenkel and Turaev seem not to have been aware of the work of
Wilson. Spiridonov and Zhedanov \cite{sz,sz2} gave an independent approach
to elliptic $6j$-symbols, showing in particular that they are
biorthogonal rational functions, and that they coincide with Wilson's
functions  in the trigonometric  limit. More precisely, 
trigonometric $6j$-symbols  correspond to certain discrete
restrictions on the
parameters of Wilson's  functions. 
Similarly, elliptic $6j$-symbols correspond to special parameter
choices for Spiridonov's and Zhedanov's biorthogonal rational
functions. We will be concerned with the larger parameter range, although,
for simplicity, we will use the term ``$6j$-symbol'' also in that setting. 

\begin{figure}[htb]
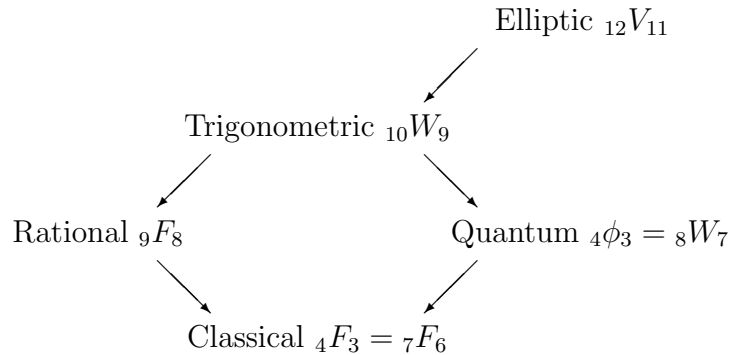
\label{figure}
\begin{diagram}[small]
 &&&&\text{\!\!\!\!\! Elliptic}\  {}_{12}V_{11}\\
&&&\ldTo &\\
&&\text{Trigonometric} \ {}_{10}W_9 &&\\
& \ldTo && \rdTo &\\
\text{Rational}\ {}_9F_8 & & & & \text{Quantum}\ {}_4\phi_3={}_8W_7\\
& \rdTo  && \ldTo & \\
&& \text{Classical}\ {}_{4}F_3={}_7F_6 &&
\end{diagram}
\caption{The hierarchy of $6j$-symbols}
\end{figure}

To summarize,  we have  a  scheme (in the sense of Askey) 
consisting of
classical, quantum, rational, trigonometric and elliptic
$6j$-symbols, see Figure \ref{figure}. Arrows  indicate 
limit transitions. We also give the
hypergeometric type of the systems. (We use Spiridonov's \cite{sp}
more  logical
notation ${}_{12}V_{11}$, see \eqref{vd} below,  
rather than ${}_{10}\omega_9$ as in \cite{ft2},   
 for the series underlying elliptic $6j$-symbols.) 

Note that the discrete part of the  Askey Scheme lies below the
classical and quantum $6j$-symbols in Figure  \ref{figure}. We remark that,
once we decide to include biorthogonal rational functions, many
further limit cases  exist (including biorthogonal polynomials and
orthogonal rational functions). It seems desirable to
classify all limit cases, along with their continuous relatives. 
Many known systems  
(see \cite{ai,av,gm,im1,im2,kw3,p,rm,rm1,rm0,rm2} for some
candidates) should fit into this larger picture.

The aim of the present work is to give a self-contained and
elementary
approach to $6j$-symbols, which works for all five cases. We will show
how to obtain many of their  properties in an elementary
fashion, without using quantum groups or techniques from statistical
mechanics (although the approach is certainly related to both).
In the exposition we will focus on  trigonometric $6j$-symbols.  We stress that
this is not because of any essential difficulties with the elliptic
case,  but since we want
to emphasize the elementary nature of our approach as much as possible.

Our main idea comes from the interpretation of Askey--Wilson
and  $q$-Racah polynomials given 
in \cite{r}; see \cite{st,z} for related work. The
standard definition of $6j$-symbols involves three-fold tensor products
of representations. This works equally well in the classical and
quantum case. In \cite{r}, we gave an interpretation of
$q$-Racah  polynomials involving a single irreducible representation
of the  $\mathrm{SL}(2)$ quantum group. On the level of polynomials, this means
that $q$-Racah polynomials appear as  $q$-analogues of Krawtchouk
polynomials rather than of Racah polynomials. Realizing the
 representation using 
 difference operators  on a function space (sometimes called
the coherent state method \cite{j}), this yields a kind of
generating function for $q$-Racah polynomials, see \eqref{red} below. 
We may now forget about the quantum
group and use the generating function to recover the main 
properties of $q$-Racah polynomials. Our aim is to generalize
this approach to include all $6j$-symbols in Figure \ref{figure}, keeping the
underlying quantum
group (known as the Sklyanin algebra in the most general case) 
implicit until the final Section~\ref{ssa}.

The plan of the paper is as follows. Section \ref{ps} contains
preliminaries; in particular we explain in some detail the degenerate cases
corresponding to Krawtchouk polynomials and $q$-Racah polynomials
(or quantum $6j$-symbols). In Section \ref{st} we generalize this to
trigonometric $6j$-symbols, and in
Section \ref{sell} we sketch the straight-forward extension 
to elliptic $6j$-symbols. 
Although the main point of the paper is to avoid using quantum groups,
we give an algebraic interpretation of our construction in
Section \ref{ssa}. It turns out that elliptic $6j$-symbols appear as
the transition matrix between the solutions of two different
generalized eigenvalue problems in a finite-dimensional representation
of the Sklyanin algebra.

{\bf Acknowledgements:} It is a pleasure to dedicate this paper to
Richard Askey, who has provided help and support in times when it was  needed.
It is a natural continuation of  \cite{r}, which was presented at
 Askey's  65th  birthday meeting in
1998, and it is directly inspired by his comments on that
talk. Furthermore, I
thank Eric Rains for crucial correspondence on the Sklyanin algebra
and Paul Terwilliger for many useful
discussions. 
The figures were made using Paul Taylor's package
\verb+diagrams+.   Finally, I 
thank Institut Mittag-Leffler for valuable support, giving
me the time to type up the manuscript.

\section{Preliminaries}
\label{ps}

\subsection{Notation}

 We  recall the standard notation for shifted factorials
$$(a)_k=a(a+1)\dotsm(a+k-1),$$
$$(a_1,\dots,a_n)_k=(a_1)_k\dotsm(a_n)_k,$$
for hypergeometric series
$${}_{r}F_{s}\!\left[\begin{matrix}a_1,\dots,a_r\\
b_1,\dots,b_s \end{matrix}\,;x\right]=
\sum_{k=0}^\infty\frac{(a_1,\dots,a_r)_k}{(b_1,\dots,b_s)_k}\frac{x^k}{k!},$$
for $q$-shifted
factorials
\begin{equation}\label{qp}(a;q)_k=(1-a)(1-aq)\dotsm(1-aq^{k-1}),
\end{equation}
$$(a_1,\dots,a_n;q)_k=(a_1;q)_k\dotsm(a_n;q)_k, $$
for $q$-binomial coefficients
$$\qb{N}{k}=\frac{(q;q)_N}{(q;q)_k(q;q)_{N-k}}, $$
for basic hypergeometric series
$${}_{r}{\phi}_s\!\left[\begin{matrix} a_1,\dots,
  a_{r}\\b_1,\dots,b_s\end{matrix};q,z\right]=
  \sum_{k=0}^{\infty}\frac{(a_1,\dots,a_{r};q)_k}
{(q,b_1,\dots,b_s;q)_k}
  \left((-1)^kq^{\binom k2}\right)^{1+s-r}z
  ^k$$
and for very-well-poised series
\begin{equation*}\begin{split}{}_{r+1}W_r(a;b_1,\dots,b_{r-2};q,z)&
  ={}_{r+1}{\phi}_r\!\left[\begin{matrix}
    a,qa^{\frac 12}, -qa^{\frac 12}, b_1,\dots,
  b_{r-2}\\a^{\frac 12},-a^{\frac
    12},aq/b_1,\dots,aq/b_{r-2}\end{matrix};q,z\right]\\
  &=\sum_{k=0}^{\infty}\frac{1-aq^{2k}}{1-a}\frac{(a,b_1,\dots,b_{r-2};q)_k}
{(q,aq/b_1,\dots,aq/b_{r-2};q)_k}\,z^k.
\end{split}\end{equation*}
If one of the numerator parameters equals $q^{-n}$, with $n$ a non-negative
integer, the series reduces to a finite sum.
We are particularly interested in the terminating balanced ${}_{10}W_9$, 
that is, the case when $r=9$, the sum is finite,
$z=q$ and $a^3q^2=b_1\dotsm b_7$.
The standard reference for all this is \cite{gr}.

To write our results in standard notation, 
some routine computation  involving $q$-shifted factorials
is necessary. We will not give any details, but
 we mention that all that one needs is the elementary identities
\begin{subequations}\label{eli}
\begin{align}
(a;q)_n&=(-1)^nq^{\binom n2}a^n(q^{1-n}/a;q)_n,\\
(a;q)_{n+k}&=(a;q)_n(aq^n;q)_k,\\
(a;q)_{n-k}&=(-1)^kq^{\binom k2}(q^{1-n}/a)^k
\frac{(a;q)_n}{(q^{1-n}/a;q)_k}. 
\end{align}
\end{subequations}

\subsection{An extended example: Krawtchouk polynomials}
\label{sk}

Our  guiding example will be Krawtchouk polynomials, arising as matrix
elements of $\mathrm{SL}(2,\mathbb C)$. (Incidentally, they  also
appear as $6j$-symbols, namely, of the oscillator algebra \cite[\S
8.6.6]{vk}.) 
For later comparison, we 
recall some fundamental facts on this topic \cite{kw,vk}.

Consider the coefficients $K_k^l=K_k^l(a,b,c,d;N)$ in
\begin{equation}\label{ke}(ax+b)^k(cx+d)^{N-k}=\sum_{l=0}^N K_k^l\,
x^l,
\end{equation}
where $k\in\{0,1,\dots,N\}$ and we assume, with no great loss of generality,
that $ad-bc=1$. Note that $\mathrm{SL}(2)$ acts on
polynomials of degree $\leq N$ by
\begin{equation}\label{gr}
p(x)\mapsto(cx+d)^Np\left(\frac{ax+b}{cx+d}\right), \end{equation}
and that $K_k^l$ are  the matrix elements of this group action in the
standard basis of monomials.

Using the binomial theorem, several different expressions for $K_k^l$ 
as hypergeometric sums may be derived. For instance,
\begin{equation*}\begin{split}
(ax+b)^k(cx+d)^{N-k}&=\left(\frac 1d\, x+\frac
bd\,(cx+d)\right)^k(cx+d)^{N-k}\\
&=\sum_{j=0}^k\binom{k}{j}\frac{b^{k-j}}{d^k}\,x^j(cx+d)^{N-j}\\
&=\sum_{j=0}^k\sum_{m=0}^{N-j}\binom{k}{j}\binom{N-j}{m}b^{k-j}c^md^{N-k-j-m}
x^{m+j}.
\end{split}
\end{equation*}
Thus, writing $m=l-j$, we obtain
\begin{equation}\label{khe}\begin{split}K_k^l&
=\sum_{j=0}^{\min(k,l)}\binom{k}{j}\binom{N-j}{l-j}b^{k-j}c^{l-j}d^{N-k-l}\\
&=\binom{N}{l}b^kc^ld^{N-k-l}\,{}_{2}F_{1}\!\left[\begin{matrix}-k,-l\\
-N \end{matrix}\,;\,-\frac{1}{bc}\right],\end{split}\end{equation}
in standard hypergeometric notation.

Note that the expansion problem inverse to \eqref{ke},
\begin{equation}\label{ike}x^k=\sum_{l=0}^N\tilde K_k^l\,
(ax+b)^l(cx+d)^{N-l}
\end{equation}
is equivalent to the original problem (replace the matrix
$\left(\begin{smallmatrix} a&b\\c&d\end{smallmatrix}\right)$
by its inverse). Thus, $\tilde K_k^l$ is given by a similar formula,
namely,
$$\tilde K_k^l=\binom{N}{l}(-1)^{k+l}a^{N-k-l}b^kc^l\,{}_{2}F_{1}\!\left[\begin{matrix}-k,-l\\
-N \end{matrix}\,;\,-\frac{1}{bc}\right]. $$
Combining \eqref{ke} and \eqref{ike}, we obtain the orthogonality relation
\begin{equation}\label{keie}\begin{split}
\delta_{km}&=\sum_{l=0}^N K_k^l\,\tilde K_l^m\\
&=\sum_{l=0}^N\binom{N}{l}\binom{N}{m}(-1)^{l+m}a^{N-m-l}b^{k+l}c^{l+m}
d^{N-k-l}\\
&\quad\times\,{}_{2}F_{1}\!\left[\begin{matrix}-k,-l\\
-N \end{matrix}\,;\,-\frac{1}{bc}\right]
\,{}_{2}F_{1}\!\left[\begin{matrix}-m,-l\\
-N \end{matrix}\,;\,-\frac{1}{bc}\right].
\end{split}\end{equation}

Now let us introduce the standard notation
$$K_n(x;p,N)={}_{2}F_{1}\!\left[\begin{matrix}-n,-x\\
-N \end{matrix}\,;\,\frac{1}{p}\right]. $$
This is a polynomial in $x$ of degree $n$, known as the
Krawtchouk  polynomial \cite{khp}.  
Writing $bc=-p$, $ad=1-p$ and $t=cx/d$, \eqref{ke} takes the form
\begin{equation}\label{kgf}
\left(1+\frac{p-1}{p}\,t\right)^k(1+t)^{N-k}=\sum_{l=0}^N
\binom{N}{l}K_l(k;p,N)\,t^l,\end{equation}
which is a well-known generating function for Krawtchouk
polynomials. Our  approach to the $6j$-symbols in Figure \ref{figure} will
be based  on generalizing this identity.

In terms of Krawtchouk polynomials, 
\eqref{keie} takes the form
$$\sum_{x=0}^N\binom{N}{x}p^x(1-p)^{N-x}K_k(x;p,N)K_m(x;p,N)=\delta_{km}
\frac{(1-p)^k}{p^k\binom{N}{k}}.$$
For $0<p<1$, this is an orthogonality relation
for a positive measure, namely, the binomial distribution on a finite 
arithmetic progression. That we get a genuine orthogonality stems from
the fact that the underlying representation is unitarizable for the group
$\mathrm{SU}(2)$.

Several other interesting properties of Krawtchouk polynomials are
immediately obtained from \eqref{ke}. For instance,
one may consider three  bases $e_k$, $f_k$, $g_k$, each being
of the form $(ax+b)^k(cx+d)^{N-k}$, with different $a$, $b$, $c$,
$d$. The transition coefficients in
$$e_k=\sum_{l}K_{kl}\,g_l=\sum_lK_{kl}'\,f_l,\qquad f_k=\sum K_{kl}''\, g_l $$
are then all given by Krawtchouk polynomials, with different parameter
$p$. Clearly, they are related by  matrix multiplication:
\begin{equation}\label{kadd}K_{nm}=\sum_{k=0}^NK_{nk}'K_{km}''.\end{equation}
In the case $e_k=g_k$, one gets back the orthogonality \eqref{keie}.

From the viewpoint of group theory, \eqref{kadd} corresponds to representing
the group law (i.e.\ multiplication of $2\times 2$ matrices) in 
an $(N+1)$-dimensional representation. 
 This should be quite familiar when
$N=1$ and we restrict to $\mathrm{SO}(2)$, the rotations of the
plane, obtaining in this way the addition formulas for sine and
cosine. Thus, \eqref{kadd} appears as a natural extension of these
addition formulas. 

From the hypergeometric viewpoint, \eqref{kadd} is an instance of
 Meixner's formula \cite{m}
\begin{multline*}
\sum_{k=0}^\infty\frac{(c)_k}{k!}
\,{}_{2}F_{1}\!\left[\begin{matrix}-k,a\\
c \end{matrix}\,;x\right]\,{}_{2}F_{1}\!\left[\begin{matrix}-k,b\\
c \end{matrix}\,;y\right] z^k\\
=\frac{(1-z)^{a+b-c}}{(1-z+xz)^a(1-z+xz)^b}
\,{}_{2}F_{1}\!\left[\begin{matrix}a,b\\
c \end{matrix}\,;\,\frac{xyz}{(1-z+xz)(1-z+yz)}\right].
\end{multline*}
More precisely, it is 
 the special case when $a=-n$, $b=-m$,  $c=-N$, with $m$, $n$, $N$
integers such that $0\leq m,n\leq N$.

Another consequence of \eqref{ke} is obtained by exploiting the
multiplicative structure of the basis vectors. Namely, expanding both
sides of
\begin{equation}\label{mb}
(ax+b)^{k+j}(cx+d)^{M+N-k-j}=(ax+b)^k(cx+d)^{M-k}(ax+b)^j(cx+d)^{N-j}
\end{equation}
into monomials gives
$$\sum_{l} K_{k+j}^l x^l=\sum_m K_k^m x^m\sum_n K_j^n x^n, $$
or
\begin{equation}\label{kconv}
K_{k+j}^l(a,b,c,d;M+N)=\sum_{m+n=l}K_k^m(a,b,c,d;M)K_j^n(a,b,c,d;N).
\end{equation}
In hypergeometric notation, this is
\begin{multline}\label{hypconv}
\binom{M+N}{l}\,{}_{2}F_{1}\!\left[\begin{matrix}-l,-k-j\\
-M-N \end{matrix}\,;t\right]\\
=\sum_{\substack{m+n=l\\0\leq m\leq M\\0\leq n\leq N}}
\binom{M}{m}\binom{N}{n}\,{}_{2}F_{1}\!\left[\begin{matrix}-m,-k\\
-M\end{matrix}\,;t\right]\,{}_{2}F_{1}\!\left[\begin{matrix}-n,-j\\
-N \end{matrix}\,;t\right].
 \end{multline}

The group-theoretic interpretation of \eqref{kconv} is the
following. Let $V_N$ denote the $(N+1)$-dimensional irreducible 
representation of $\mathrm{SL}(2)$, realized on the space of
polynomials as above. Then multiplication of polynomials defines a
map $V_M\otimes V_N\rightarrow V_{M+N}$. The relation
\eqref{mb}, and thus
\eqref{kconv}, expresses the fact that this map is intertwining, that
is, commutes with the group action. This  immediately
suggests a non-trivial generalization. Namely, one has the equivalence
of representations $V_M\otimes
V_N\simeq\bigoplus_{s=0}^{\min(M,N)}V_{M+N-2s}$, and one may do the
same thing for the intertwiners  $V_M\otimes V_N\rightarrow V_{M+N-2s}$.
The corresponding generalization of \eqref{hypconv} has additional
factors of type ${}_3F_2$ appearing on both sides. 
From the group-theoretic viewpoint, these  are Clebsch--Gordan
coefficients and, from the viewpoint  of special functions, Hahn
polynomials, see \cite[\S 8.5.3]{vk}.

Next we point out that \eqref{kconv} may be iterated to
\begin{multline}\label{mkc}
K_{k_1+\dots+k_n}^l(a,b,c,d;M_1+\dots+M_n)\\
=\sum_{m_1+\dots+m_n=l}K_{k_1}^{m_1}(a,b,c,d;M_1)\dotsm
K_{k_n}^{m_n}(a,b,c,d;M_n),
\end{multline}
where $0\leq k_i, m_i\leq M_i$.
This is especially interesting when $M_i=1$  for
all $i$. Writing the result in
hypergeometric form, we get in that case
$$
\binom{n}{l}\,{}_{2}F_{1}\!\left[\begin{matrix}-l,-k_1-\dots-k_n\\
-n \end{matrix}\,;t\right]
=\sum_{\substack{m_1+\dots+m_n=l\\0\leq m_i\leq 1}}
\prod_{i=1}^n\,{}_{2}F_{1}\!\left[\begin{matrix}-m_i,-k_i\\
-1\end{matrix}\,;t\right].
$$
Note that the range of summation may be identified with the
$l$-element subsets $L$ of $N=\{1,\dots,n\}$  (interpreting $m_i=1$
as $i\in L$). Similarly, $(k_1,\dots,k_n)$ labels a subset $K$ of 
$N$ with $\sum_i k_i$ elements. 
Since
$$\,{}_{2}F_{1}\!\left[\begin{matrix}-m_i,-k_i\\
-1\end{matrix}\,;t\right]=\begin{cases}
1-t, & m_i=k_i=1,\\
1, & \text{otherwise},
\end{cases}
$$
the term in the sum is  $(1-t)^{|L\cap K|}$. Replacing $t$ with
$1-t$,   we now have 
\begin{equation}\label{kcomb}
\binom{|N|}{l}\,{}_{2}F_{1}\!\left[\begin{matrix}-l,-|K|\\
-|N| \end{matrix}\,;1-t\right]
=\sum_{L\subseteq N,\ |L|=l}t^{|L\cap K|},\qquad K\subseteq N. 
  \end{equation}
This (not very deep)
 identity  gives a combinatorial
interpretation for Krawtchouk polynomials as a generating function for
 the statistics $|L\cap K|$ on  subsets $L$ of fixed cardinality.
We shall see that the appearance of $6j$-symbols in statistical
 mechanics is via a generalization of this identity.

\subsection{$q$-Racah polynomials}
\label{rss}

In \cite{r}, we
considered a $q$-analogue of the above set-up, leading to general $q$-Racah
polynomials. The
group $\mathrm{SL}(2)$ was replaced by a quantum group, and the basis
vectors $x^k$ and $(ax+b)^k(cx+d)^{N-k}$ by appropriate $q$-shifted
products such as 
$$\prod_{j=0}^k(axq^j+b)\prod_{j=0}^{N-k}(cxq^j+d). $$
Such bases 
were interpreted as eigenvectors of  Koornwinder's 
twisted primitive elements \cite{kw2}, and also as the image of
standard basis vectors $x^k$ under Babelon's vertex-IRF transformations
\cite{b} (called generalized group elements in \cite{r}). Actually,
we focused on the case of infinite-dimensional
representations,  and only
mentioned the case of present interest somewhat parenthetically
\cite[Section 6]{r}. 

To be more precise, the expansion problem that yields quantum
$6j$-symbols  ($q$-Racah polynomials) is
\begin{equation}\label{re}\prod_{j=0}^k(axq^{-j}+b)\prod_{j=0}^{N-k}
(cxq^{-j}+d)=
\sum_{l=0}^N C_k^l\, \prod_{j=0}^l(\al xq^{j}+\be)\prod_{j=0}^{N-l}(\ga
xq^{j}+\de).\end{equation}
For generic parameter values, the polynomials on the right form a
basis for the space of polynomials of degree $\leq N$, so that the
coefficients exist uniquely. For the rest of this section we assume
that we are in such a generic situation.

Note that when $q=1$, \eqref{re} reduces to
$$(ax+b)^k(cx+d)^{N-k}=
\sum_{l=0}^N C_k^l\,(\al x+\be)^l(\ga x+\de)^{N-l},$$
which is further reduced to \eqref{ke} by a change of variables. The
expansion \eqref{re} is more rigid. After multiplying with a trivial factor and
changing parameters, we may restrict to the case
\begin{equation}\label{red}(ax;q^{-1})_k(bx;q^{-1})_{N-k}=
\sum_{l=0}^N C_k^l(a,b,c,d;N;q) \, (cx;q)_l(dx;q)_{N-l}.
\end{equation}
We could dilate $x$ to get rid of one more parameter,
 but the remaining $7$  parameters, counting $q$, enter in a
non-trivial fashion. Indeed, we have
\begin{equation}\label{cex}\begin{split}
C_k^l(a,b,c,d;N;q)&=q^{l(l-N)}\qb{N}{l}\frac{(q^{1-N}b/d;q)_l
 (q^{1-N}b/c;q)_{N-l} (q^{1-k}a/c;q)_k}{(q^{l-N}c/d;q)_l
 (q^{-l}d/c;q)_{N-l} (q^{1-N}b/c;q)_k}\\ 
&\quad\times {}_{4}{\phi}_3\!\left[\begin{matrix} q^{-k},q^{-l},q^{k-N}b/a,
 q^{l-N}c/d\\q^{-N},c/a,q^{1-N}b/d\end{matrix};q,q\right].
\end{split}\end{equation}
In Section \ref{st} we will  derive a more general identity in an
 elementary way.

Similarly as for \eqref{ke}, we may invert
\eqref{red} to get the orthogonality relation
\begin{equation}\label{ro}\delta_{km}=
\sum_{l=0}^N C_k^l(a,b,c,d;N;q)\,C_l^m(c,d,a,b;N;q^{-1}).
\end{equation}
One may verify that \eqref{ro} gives the orthogonality of
$q$-Racah polynomials. (If we want a positive
measure, some conditions on the parameters must be imposed.)

Note that \eqref{red} generalizes the generating function \eqref{kgf}
to the level of $q$-Racah
polynomials. This identity was  obtained, in
a  related but not
identical context, by Koelink and Van der Jeugt \cite[Remark 4.11(iii)]{kv}.

The mixture of base $q$ and $q^{-1}$ in \eqref{red} is
crucial. Admittedly, the expansion problem
\begin{equation}\label{rfe}(ax;q)_k(bx;q)_{N-k}=
\sum_{l=0}^N D_k^l(a,b,c,d;N;q) \, (cx;q)_l(dx;q)_{N-l}
\end{equation}
is immediately reduced to \eqref{red} by a change of variables;
explicitly, one has
$$D_k^l(a,b,c,d;N;q)=C_k^l(aq^{k-1},bq^{N-k-1},c,d;N;q).$$
However, the relation
$$\delta_{km}=
\sum_{l=0}^k D_k^l(a,b,c,d;N;q)\,D_l^m(c,d,a,b;N;q)
$$
is not equivalent to \eqref{ro}. It gives a system of biorthogonal
rational functions. When $c\bar d=a\bar b\in\mathbb R$
 it is an orthogonal system, 
found in an equivalent context by Koelink \cite[Proposition 9.5]{k};
see   also \cite[Corollary 4.4]{gm}  
and Remark \ref{dgr} below.

\begin{remark}
We conclude the introductory part of the paper with some comments on
the relation to Terwilliger's concept of a Leonard pair; see \cite{to}
and references given there. As was mentioned above, if we let  
 $e_k=(ax;q^{-1})_k(bx;q^{-1})_{N-k}$, $f_k=(cx;q)_k(dx;q)_{N-k}$,
then $e_k$ and $f_k$ appear as eigenbases of certain $q$-difference
operators $Y_1$, $Y_2$, respectively. It is easy to check that each of
these operators acts tridiagonally on the eigenbasis of the other,
that is, 
$$Y_1f_k\in\operatorname{span}\{f_{k-1},f_k,f_{k+1}\},\qquad
Y_2e_k\in\operatorname{span}\{e_{k-1},e_k,e_{k+1}\}.$$
Except for a non-degeneracy condition, this is the definition of
$(Y_1,Y_2)$ being a Leonard pair. Then \eqref{cex} means that 
$(Y_1,Y_2)$
 is a Leonard pair of ``$q$-Racah type'', which is the most general kind. 
This gives a simple model for studying Leonard pairs. For instance, 
the ``split decompositions'' \cite{ts} are easily
understood in this model. A typical split basis between $e_k$ and
$f_k$ would be $g_k=(ax;q^{-1})_k(dx;q)_{N-k}$, which interpolates
between the two other bases in the sense that
$$g_k\in\operatorname{span}\{e_k,e_{k+1},\dots,e_N\}\cap\operatorname{span}
\{f_0,f_1,\dots,f_k\}.$$
The corresponding split decomposition is then simply
$\bigoplus_{k=0}^N\mathbb C g_k$. More generally, we may picture the
 factors $(ax;q^{-1})_k$, $(bx;q^{-1})_k$, $(cx;q)_{N-k}$,
$(dx;q)_{N-k}$ as being attached to the corners of a tetrahedron, with
two edges corresponding to the original Leonard pair and four edges
corresponding to different split decompositions.
\end{remark}

\section{Trigonometric $6j$-symbols}
\label{st}

It is not hard to check that both $q$-Racah polynomials and Koelink's
orthogonal functions are  degenerate cases of
Wilson's biorthogonal functions. 
Thus, if one wants to obtain  general trigonometric  $6j$-symbols
in a similar way, it seems necessary to unify the products $(ax;q)_k$
and $(ax;q^{-1})_k$. 
The correct unification turns out to be the Askey--Wilson monomials
$h_k(x;a)=h_k(x;a;q)$, which are the natural building blocks of 
 Askey--Wilson polynomials \cite{aw2}. They are given by
$$h_k(x;a)=\prod_{j=0}^{k-1}(1-axq^j+a^2q^{2j}).$$
To write this in the notation \eqref{qp} one must introduce an auxiliary
variable $\xi$ satisfying 
\begin{equation}\label{xi}\xi+\xi^{-1}=x;\end{equation}
 then
\begin{equation}\label{awm}h_k(x;a)=(a\xi,a\xi^{-1};q)_k.\end{equation}
(In the context of Askey--Wilson polynomials one usually dilates $x$
by a factor $2$ and writes $x/2=\cos\theta$, $\xi=e^{i\theta}$.)
We will need the elementary identities
\begin{equation}\label{hi}h_k(x;a)=q^{k(k-1)}a^{2k}h_k(x;q^{1-k}/a),
\end{equation}
\begin{equation}\label{hpf}h_{k+l}(x;a)=h_k(x;a)h_l(x;aq^k).\end{equation}

It is easy to see that
\begin{equation}\label{lim}\lim_{t\rightarrow 0} h_k(x/t;at)=(ax;q)_k,\qquad
\lim_{t\rightarrow 0}
t^{2k}h_k(x/t;a/t)=a^{2k}q^{k(k-1)}(x/a;q^{-1})_k.\end{equation}
Thus, we may unify \eqref{red} and \eqref{rfe}, together with several
related expansion problems (see Remark \ref{dgr} below), into 
\begin{equation}\label{tre}
h_k(x;a)h_{N-k}(x;b)=\sum_{l=0}^N R_k^l(a,b,c,d;N;q)\,
 h_l(x;c)h_{N-l}(x;d).\end{equation}
We will suppress parameters when convenient, writing
$$R_k^l=R_k^l(a,b,c,d;N)=R_k^l(a,b,c,d;N;q). $$
Note that, in contrast to the limit cases considered above, we cannot
get rid of any parameters by scaling $x$. We shall see that $R_k^l$
depends on all $8$ parameters (counting $q$) in a non-trivial fashion.

Clearly, the coefficients $R_k^l$ exist uniquely if and only if
$(h_k(x;c)h_{N-k}(x;d))_{k=0}^N$ form a basis for the space of
polynomials of degree $\leq N$. Although it is not quite necessary for
our purposes (see Remark \ref{pbr}), we will first settle this question.

\begin{lemma}\label{bl}
The polynomials $(h_k(x;c)h_{N-k}(x;d))_{k=0}^N$ form a basis for
the space of 
polynomials of degree at most $N$ if and only if none of the following
conditions are satisfied:
\begin{subequations}\label{c}
\begin{equation}\label{c1}
c/d\in\{q^{1-N},q^{2-N},\dots,q^{N-1}\},  \end{equation}
\begin{equation}\label{c2}cd\in\{q^{1-N},q^{2-N},\dots,1\},  \end{equation}
\begin{equation}\label{c3}c=d=0. \end{equation}
\end{subequations}
\end{lemma}

\begin{proof}
If $c/d=q^j$ with $1-N\leq j\leq 0$,
then all the polynomials have the
common zero $x=d+d^{-1}$, so they cannot form a basis. Similarly, if
$c/d=q^j$ with $0\leq j\leq N-1$ 
then $x=c+c^{-1}$ is a common zero, and if \eqref{c2} holds
then both $x=c+c^{-1}$ and $x=d+d^{-1}$ are common zeroes.
In the case \eqref{c3}, all the
polynomials  equal  $1$ and clearly do not form a basis.

Conversely, assume that  none of the conditions \eqref{c} hold.
We need to show that any  linear relation
\begin{equation}\label{lr}
\sum_{k=0}^N \lambda_k\, h_k(x;c)h_{N-k}(x;d)\equiv 0\end{equation}
is trivial. By symmetry, we may
assume $c\neq 0$. 
Choosing $x=c+c^{-1}$ in \eqref{lr} gives
$$\lambda_0(dc,d/c;q)_N=0. $$
Since $(dc,d/c;q)_N=0$ only if \eqref{c1} or \eqref{c2} holds, we have
 $\lambda_0=0$. We may then divide \eqref{lr} with $1-cx+c^2$,
giving
$$\sum_{k=1}^N \lambda_k\, h_{k-1}(x;cq)h_{N-k}(x;d)\equiv 0.$$
By iteration (choosing $x=cq+(cq)^{-1}$ in the next step) or by
induction of $N$, we 
conclude that $\lambda_i=0$ for all $i$, and thus that
the polynomials form a basis.
\end{proof}

We  now turn to the problem of computing the coefficients 
$R_k^l$.
Recall that our derivation of \eqref{khe} consisted in 
applying the binomial theorem twice.  The 
same proof should be applicable to \eqref{tre}, once we have a
generalized binomial theorem of the form
\begin{equation}\label{gbt}
h_N(x;a)=\sum_{k=0}^N C_k^N(a,b,c)\, h_k(x;b)h_{N-k}(x;c).\end{equation}
In fact, \eqref{gbt} is solved by  one of the most
fundamental results on basic hypergeometric series: Jackson's
${}_8W_7$ summation \cite{gr,ja}. Since we have promised to give a
self-contained treatment, we give a straight-forward proof, 
motivated by our present view of \eqref{gbt} as an extension of
the binomial theorem. 

We will follow the standard
inductive proof of the binomial theorem based on Pascal's triangle. 
First we write
$$h_{N+1}(x;a)=h_N(x;a)(1-aq^Nx+a^2q^{2N}).$$
To get a recurrence for  $C_k^N$, we must split the
factor $1-axq^N+a^2q^{2N}$ into parts that attach to the right-hand
side of \eqref{gbt}, that is, as
\begin{equation}\label{lfs}
1-aq^Nx+a^2q^{2N}=A_k(1-bq^kx+b^2q^{2k})+B_k(1-cq^{N-k}x+c^2q^{2(N-k)}).
\end{equation}
We compute
\begin{equation}\label{rc}
\begin{split}
A_k&=\frac{(1-acq^{2N-k})(1-aq^k/c)}{(1-bcq^N)(1-bq^{2k-N}/c)},\\
B_k&=\frac{(1-abq^{N+k})(1-aq^{N-k}/b)}{(1-bcq^N)(1-cq^{N-2k}/b)},
\end{split}\end{equation}
assuming that the denominators are non-zero. For the elliptic extension 
discussed in Section \ref{sell} it is important to note that this uses
the elementary identity
\begin{multline}\label{tadd}
\frac vx\,(1-xy)(1-x/y)(1-uv)(1-u/v)\\
=(1-ux)(1-u/x)(1-vy)(1-v/y)-(1-uy)(1-u/y)(1-vx)(1-v/x),
\end{multline}
with
$$(u,v,x,y)\mapsto(cq^{N-k},bq^k,aq^N,\xi).$$

Combining \eqref{gbt} and \eqref{lfs}  yields the generalized Pascal triangle
\begin{equation}\label{gpt}C_k^{N+1}=B_kC_k^N+A_{k-1}C_{k-1}^N,\end{equation}
with boundary conditions
$$C_0^0=1,\quad C_{-1}^N=C_{N+1}^N=0. $$
Iterating \eqref{gpt}, one  quickly  guesses that
\begin{equation}\label{gbc}
C_k^N=q^{k(k-N)}\qb{N}{k}\frac{(a/c,q^{N-k}ac;q)_k(a/b,q^kab;q)_{N-k}}
{(q^{k-N}b/c;q)_k(q^{-k}c/b;q)_{N-k}(bc;q)_N}.\end{equation}
To verify the guess, we plug \eqref{gbc} into \eqref{gpt}. After
 cancelling  common factors, we are left with 
\begin{multline*}
q^{k-N-1}(1-q^N)(1-q^Nab)(1-q^Nac)(1-q^{N+1-2k}c/b)\\
\begin{split}&=(1-q^k)(1-q^{k-1}ab)(1-q^{2N-k+1}ac)(1-q^{-k}c/b)\\
&\quad-(1-q^{k-N-1})(1-q^{N+k}ab)(1-q^{N-k}ac)(1-q^{N+1-k}c/b),
\end{split}\end{multline*}
which is another instance of \eqref{tadd}, this time with
$$(u,v,x,y)\mapsto(q^{N+\frac 12}\sqrt{ac},q^{-\frac
12}\sqrt{ac},q^{N-k+\frac 12}\sqrt{ac},q^{k-\frac 12}b\sqrt{a/c}).$$
This shows that,  for generic parameters,
 \eqref{gbt} holds with the coefficients given by \eqref{gbc}.

\begin{remark}
\label{pbr}
Note that our proof 
 did not use Lemma \ref{bl}. We see from the computation
 that the expansion exists uniquely unless there are
zeroes in the denominators of \eqref{rc}, which happens precisely if
$$b/c\in\{q^{1-N},q^{2-N},\dots,q^{N-1}\},\qquad 
bc\in\{1,q^{-1},\dots,q^{1-N}\} \qquad \text{or}\quad b=c=0.$$
As expected, this corresponds exactly to the conditions \eqref{c}.
\end{remark}

\begin{remark}
Plugging \eqref{gbc} into \eqref{gbt} and rewriting the result in
standard  notation gives
$${}_8W_7(q^{-N}b/c;q^{-N},q^{1-N}/ac,a/c,b\xi,b\xi^{-1};q,q)= 
\frac{(cb,c/b,a\xi,a\xi^{-1};q)_N}{(ab,a/b,c\xi,c\xi^{-1};q)_N}.$$ 
This is Jackson's summation.
Essentially the same method was used in \cite{r2} to obtain extensions
of Jackson's summation to multiple elliptic hypergeometric series
related to the root systems $A_n$ and $D_n$.
\end{remark}

We may now  compute the coefficients $R_k^l$ in \eqref{tre} by applying the
``binomial theorem'' \eqref{gbt}
twice.  For instance, using
\eqref{hpf} we may write
\begin{multline*}h_k(x;a)h_{N-k}(x;b)
=\sum_{j=0}^k C^k_j(a,c,bq^{N-k})\,
h_j(x;c)h_{N-j}(x;b)\\
=\sum_{j=0}^k\sum_{m=0}^{N-j}
C^k_j(a,c,bq^{N-k}) C_{m}^{N-j}(b,cq^j,d)\,h_{j+m}(x;c)h_{N-j-m}(x;d).
\end{multline*}
This gives
 $$R_k^l=\sum_{j=0}^{\min(k,l)}C_j^k(a,c,bq^{N-k})C_{l-j}^{N-j}(b,cq^j,d). $$
Plugging in the expressions from
\eqref{gbc} and rewriting the result in standard form
one finds that the sum is a balanced ${}_{10}W_9$ series.

\begin{theorem}
\label{ret}
 For generic values of the parameters, the 
coefficients $R_k^l$ in \eqref{tre} exist uniquely and  are given by
\begin{multline*}
R_k^l(a,b,c,d;N;q)\\
\begin{split}&=q^{l(l-N)}\qb{N}{l}\frac{(ac,a/c;q)_k(q^{N-l}bd,b/d;q)_l
(b/c;q)_{N-k}(b/c;q)_{N-l}
(bc;q)_{N-k}}{(q^{l-N}c/d;q)_l(q^{-l}d/c;q)_{N-l}(cd;q)_N(b/c;q)_N(bc;q)_l}\\
&\quad
\times{}_{10}W_9(q^{-N}c/b;q^{-k},q^{-l},q^{k-N}a/b,q^{l-N}c/d,cd,q^{1-N}/ab,
qc/b;q,q).
\end{split}\end{multline*}
\end{theorem}

\begin{remark}
The special case $d=0$ of  Theorem \ref{ret} was recently obtained by
Ismail and Stanton  \cite[Theorem 3.1]{is} using different methods.
\end{remark}

\begin{remark}
From their definition, it is clear that $R_k^l$ have the  symmetries
\begin{subequations}\label{symm}
\begin{equation}
R_k^l(a,b,c,d;N)=R_{N-k}^l(b,a,c,d;N)=R_{k}^{N-l}(a,b,d,c;N), 
\end{equation}
and from \eqref{hi} we have moreover that
\begin{equation}
R_k^l(a,b,c,d;N)=q^{-2\binom k 2}a^{-2k}R_k^l(q^{1-k}/a,b,c,d;N).\end{equation}
\end{subequations}
Combining these
symmetries with Theorem~\ref{ret} gives
further expressions for $R_k^l$ as  ${}_{10}W_9$ sums. These
 are related via Bailey's classical ${}_{10}W_9$
transformations \cite{gr}. On the other hand, the explicit expression in
Theorem \ref{ret} implies many
 symmetries for $R_k^l$ that are not obvious from the definition.  
\end{remark}

\section{Elementary properties}

\subsection{Biorthogonality}\label{bs}

It is clear from \eqref{tre} that the coefficients $R_k^l$  satisfy 
\begin{equation}\label{tbo}
\delta_{nm}=\sum_{k=0}^N R_n^k(a,b,c,d;N;q)\,R_k^m(c,d,a,b;N;q).
\end{equation}
We will now show that \eqref{tbo} gives a system of biorthogonal rational
functions, which is identical to the one obtained by Wilson \cite{w2}.

To facilitate comparison with Wilson's result, we rewrite 
\eqref{tbo} in terms of the functions
\begin{equation*}\begin{split}R_n(\mu(k))&=\frac{q^{k(N-k)}}{(cd)^n\qb{N}{k}}
\frac{(q^{-N};q)_n(q^{k-N}c/d,bc;q)_k(q^{-k}d/c,bd;q)_{N-k}(cd;q)_N}
{(b/d;q)_k(bc,bd;q)_{N-n}(b/c;q)_{N-k}}\\
&\quad\times R_n^k(a,b,c,d;N;q)\\
&=\frac{(q^{-N},ac,q^{1-N}/bd,a/c;q)_n}{(q^{1-N}c/b;q)_n}\\
&\quad\times
{}_{10}W_9(q^{-N}c/b;q^{-n},q^{n-N}a/b,q^{-k},q^{k-N}c/d,cd,q^{1-N}/ab,
cq/b;q,q)
\end{split}\end{equation*}
and
\begin{equation*}\begin{split}S_m(\mu(k))&=\frac{q^{m(N-m)}}{(ab)^m\qb{N}{m}}
\frac{(q^{-N},ac,ad,q^{m-N}a/b;q)_m(q^{-m}b/a;q)_{N-m}(ab;q)_N}
{(ac,c/a;q)_k(ad,d/a;q)_{N-k}}\\
&\quad\times R_k^m(c,d,a,b;N;q)\\
&=\frac{(q^{-N},ac,q^{1-N}/bd,d/b;q)_m}{(q^{1-N}a/d;q)_m}\\
&\quad\times
{}_{10}W_9(q^{-N}a/d;q^{-m},q^{m-N}a/b,q^{-k},q^{k-N}c/d,ab,q^{1-N}/cd,
aq/d;q,q),
\end{split}\end{equation*}
where
$$\mu(k)=q^{-k}+q^{k-N}c/d. $$

Note that $R_n$ has the form
\begin{equation*}\begin{split}
R_n(\mu(k))&=\sum_{j=0}^{n}\sigma_j\frac{(q^{-k},q^{k-N}c/d;q)_j}
{(q^{1-N+k}c/b,q^{1-k}d/b;q)_j}\\
&=\sum_{j=0}^{n}\sigma_j
\prod_{t=0}^{j-1}\frac{1-q^t\mu(k)+q^{2t-N}c/d}{1-q^{t+1}\mu(k)d/b+
q^{2t+2-N}cd/b^2},\end{split}\end{equation*}
with $\sigma_j$ independent of $k$, and is thus  a rational
function in $\mu(k)$ of degree $n/n$. 
Similarly, 
$$S_m(\mu(k))=\sum_{j=0}^{m}\tau_j
\prod_{t=0}^{j-1}\frac{1-q^t\mu(k)+q^{2t-N}c/d}
{1-q^{t+1}\mu(k)a/c+q^{2t+2-N}a^2/cd}. $$

In terms of these functions, \eqref{tbo} takes the form
\begin{equation}\label{eto}
\sum_{k=0}^N w_k\, R_n(\mu(k))S_m(\mu(k))=C_n\,\delta_{nm},\end{equation}
where 
$$w_k=\frac{1-q^{2k-N}c/d}{1-q^{-N}c/d}
\frac{(q^{-N}c/d,q^{-N},ac,q^{1-N}/bd,b/d,c/a;q)_k}
{(q,qc/d,q^{1-N}/ad,bc,q^{1-N}c/b,q^{1-N}a/d;q)_k}\,q^k $$
and
\begin{equation*}\begin{split}
C_n&=\frac{(ba,b/a,dc,d/c;q)_N}{(bc,b/c,da,d/a;q)_N}\\
&\quad\times\frac{1-q^{-N}a/b}{1-q^{2n-N}a/b}
\frac{(q,q^{-N},ac,ad,q^{1-N}/bc,q^{1-N}/bd,aq/b;q)_n}{(q^{-N}a/b;q)_n}
\,q^{-n}.
\end{split}\end{equation*}
Thus, we have indeed a system of biorthogonal rational
functions.

We now compare this result with the work of Wilson \cite{w2}, who
 used the notation
\begin{multline}\label{wof}
r_n\left(\frac{z+z^{-1}}2;a,b,c,d,e,f;q\right)\\
=\frac{(ab,ac,ad,1/af;q)_n}{(aq/e;q)_n}
\,{}_{10}W_9(a/e;az,a/z,q/be,q/ce,q/de,q^n/ef,q^{-n};q,q),
\end{multline}
where
$$abcdef=q. $$
The normalization is chosen so as
to make $r_n$ symmetric in $a$, $b$, $c$, $d$.
Assuming $ab=q^{-N}$ with $N$ a non-negative integer, Wilson obtained the
biorthogonality relation 
\begin{multline}\label{wo}
\sum_{k=0}^N w_k\, r_n\left(\frac{aq^k+a^{-1}q^{-k}}2;a,b,c,d,e,f;q\right)\\
\times r_m\left(\frac{aq^k+a^{-1}q^{-k}}2;a,b,c,d,f,e;q\right)
=C_n\,\delta_{nm},\end{multline}
where
$$ w_k=\frac{1-a^2q^{2k}}{1-a^2}\frac{(a^2,ab,ac,ad,ae,af;q)_k}
{(q,aq/b,aq/c,aq/d,aq/e,aq/f;q)_k}\,q^k$$
and 
$$C_n=\frac{(a^2q,q/cd,q/ce,q/de;q)_N}{(aq/c,aq/d,aq/e,bf;q)_N}
\frac{(q,q^n/ef,ab,ac,ad,bc,bd,cd;q)_n}{(q/ef;q)_{2n}}\,q^{-n} $$
(in \cite{w2}, the factor $q^{-n}$ and the exponent $2$ in $a^2q$ are
missing in the expression for $C_n$).
Note that the case $m=n=0$ of \eqref{wo} is the Jackson sum.

It is now easy to check that \eqref{eto} and \eqref{wo} are
equivalent. The explicit correspondence of parameters is
\begin{multline*}(a,b,c,d,e,f)\\
\mapsto(q^{-\frac N2}\sqrt{c/d},q^{-\frac N2}\sqrt{d/c},
 q^{\frac N2}a\sqrt{cd},
q^{1-\frac N2}/b\sqrt{cd}, q^{\frac N2}b/\sqrt{cd},q^{\frac
 N2}\sqrt{cd}/a)
\end{multline*} 
(which is consistent with the  relations $ab=q^{-N}$, $abcdef=q$)
or, conversely,
\begin{equation}\label{cov}
(a,b,c,d)\mapsto(\sqrt{c/f},q/d\sqrt{cf},a\sqrt{cf},b\sqrt{cf}).
\end{equation}

\begin{remark}
Continuous biorthogonality measures for the function $r_n$ (not
assuming $ab=q^{-N}$) were obtained
by Rahman \cite{rm1,rm2}, see \cite{sz3} for the elliptic case.
\end{remark}

\begin{remark}\label{dgr}
Note that, in view of the limit relations \eqref{lim}, any one of the
sixteen expansion problems
$$(ax;q^\pm)_k(bx;q^\pm)_{N-k}=
\sum_{l=0}^N C_k^l \, (cx;q^\pm)_l(dx;q^\pm)_{N-l}, $$
with all possible choices of $\pm$, may be obtained as a degenerate
case of \eqref{tre}. It is easy to see from Theorem \ref{ret} that the
coefficients $C_k^l$ are always given by ${}_4\phi_3$ or
(equivalently, in view of Watson's transformation \cite{gr}) ${}_8W_7$
sums. Gupta and Masson \cite{gm} worked out all such degenerate cases
of Wilson's biorthogonal rational functions, finding five different systems.
 The system in \cite[Corollary 4.2]{gm} is related to the expansion
$$(ax;q)_k(bx;q^{-1})_{N-k}=
\sum_{l=0}^N C_k^l \, (cx;q)_l(dx;q^{-1})_{N-l}, $$
the system in \cite[Corollary 4.3]{gm} to 
$$(ax;q)_k(bx;q^{-1})_{N-k}=
\sum_{l=0}^N C_k^l \, (cx;q)_l(dx;q)_{N-l}, $$
the system in \cite[Corollary 4.4]{gm} is essentially Koelink's
functions \cite{k}, related to
$$(ax;q)_k(bx;q)_{N-k}=
\sum_{l=0}^N C_k^l \, (cx;q)_l(dx;q)_{N-l}, $$
the system in \cite[Corollary 4.5]{gm} to
$$(ax;q)_k(bx;q)_{N-k}=
\sum_{l=0}^N C_k^l \, (cx;q)_l(dx;q^{-1})_{N-l}, $$
and the system in \cite[Corollary 4.6]{gm} is the $q$-Racah
polynomials, related to
$$(ax;q^{-1})_k(bx;q^{-1})_{N-k}=
\sum_{l=0}^N C_k^l \, (cx;q)_l(dx;q)_{N-l}. $$
(For some of the systems Gupta and Masson gave a more
general version, with  infinite discrete biorthogonality measure.) 
All other cases may be reduced to one of those five.
\end{remark}

\subsection{Addition formula}

By iterating \eqref{tre}, one immediately generalizes the
biorthogonality relation \eqref{tbo} to
\begin{equation}\label{radd}
R_n^m(a,b,e,f;N;q)=\sum_{k=0}^N R_n^k(a,b,c,d;N;q)\,R_k^m(c,d,e,f;N;q).
\end{equation}
This is an extension  of the addition formula 
\eqref{kadd}.
 We do not believe that the general case of \eqref{radd} can be found
in the literature, although it can probably be obtained by analytic
continuation from the Yang--Baxter equation for trigono\-metric
$6j$-symbols \cite{d2,ft1}. 
Though in the present approach it seems almost trivial,
in a more direct approach, such as defining $R_k^m$ through the
explicit expression in Theorem \ref{ret}, it might not be easy to
guess the existence of such an identity, nor to give a proof.

It may be of interest to rewrite \eqref{radd} in Wilson's notation
\eqref{wof}. 
We introduce $s=e/a$, $t=b/f$ as new parameters, and then make the
change of variables \eqref{cov}.  
 The calculations are essentially
the same as those in Section \ref{bs}, and we are content with stating the
end result.

\begin{corollary}
For $abcdef=q$, $ab=q^{-N}$ and $s$ and $t$ arbitrary,
 Wilson's functions \eqref{wof} satisfy the addition formula
\begin{multline}\label{wadd}
\sum_{k=0}^N w_k\, r_n\left(\frac{aq^k+a^{-1}q^{-k}}2;a,b,c,d,e,f;q\right)\\
\begin{split}&\quad
\times r_m\left(\frac{aq^k+a^{-1}q^{-k}}2;a,b,cs,dt,f/s,e/t;q\right)\\
&=X\,R_n^m(\sqrt{c/f},q/d\sqrt{cf},s\sqrt{c/f},
q/td\sqrt{cf};q,N)\\
&=Y\, r_n\left(\frac{Aq^m+A^{-1}q^{-m}}2;A,B,C,D,E,F;q\right),
\end{split}\end{multline}
where
\begin{equation*}\begin{split}
w_k&=\frac{1-a^2q^{2k}}{1-a^2}\frac{(a^2,ab,acs,ad,ae,af/s;q)_k}
{(q,aq/b,aq/cs,aq/d,aq/e,aqs/f;q)_k}\,q^k,\\
X&=
\frac{(a^2q,q/cdst,qt/ces,q/de;q)_N}{(aq/cs,aq/d,aq/e,bf/s;q)_N}\,
{(ab,ad,bd;q)_n}\\
&\quad\times\frac{(q,q^m st/ef,acs,bcs,cdst;q)_m}{(qst/ef;q)_{2m}}\,q^{-n}
t^{2m-N}q^{m^2-n^2}(ce^2f)^{n-m},\\
Y&=\frac{(a^2q,q/de,q/cds,q/ces;q)_N}{(aq/cs,aq/d,aq/e,bf/s;q)_N}
\frac{(ad,bd;q)_n}{(q/ces,dt/e;q)_n}\\
&\quad\times\frac{(ab,t,acs,bcs,dt/e;q)_m}{(qs/df,qs/ef;q)_m},
\end{split}\end{equation*}
$$(A,B,C,D,E,F)=
(\sqrt{st/ef},ab\sqrt{ef/st},c\sqrt{es/ft},d\sqrt{ft/es},\sqrt{eft/s},
\sqrt{efs/t}).$$
\end{corollary} 

The intermediate expression in \eqref{wadd} 
makes it clear that
the special case $s=t=1$ gives back \eqref{wo}, since then
$R_n^m=\delta_{nm}$.  The presence of square
roots is due to Wilson's 
 choice of parametrization. Writing
the identity explicitly in terms of ${}_{10}W_9$-series, all square
roots combine or cancel.

\subsection{Convolution formulas}

Next we extend \eqref{kconv} to the present setting, by exploiting the
multiplicative property \eqref{hpf} of our basis elements. Because of 
the shifts appearing in that identity there are  several
different convolution formulas, which we  write compactly as follows.

\begin{corollary}\label{coco}
The coefficients $R_k^m$ satisfy the convolution formulas
\begin{multline}\label{rconv}
R_{k+j}^l(a,b,c,d;M+N;q)
=\sum_{m+n=l}R_k^m(aq^{\alpha j},bq^{\beta(N-j)},c,
d;M;q)\\
\times R_j^n(aq^{(1-\alpha)k},bq^{(1-\beta)(M-k)},
cq^{m},dq^{M-m};N;q)
\end{multline} 
for all $\alpha,\beta\in\{0,1\}$, where $0\leq k,m\leq
M$, $0\leq j,n\leq N$. 
\end{corollary}

\begin{proof}
Since, for generic parameters, $R_{k}^l$ is determined by \eqref{tre},
it suffices to compute 
$$\sum_{l=0}^{M+N} C_l\,h_l(x;c)h_{M+N-l}(x;d), $$
where $C_l$ is the right-hand side of \eqref{rconv}. 
Inside the summation sign, we split the factors as
$$ h_{m+n}(x;c)h_{M+N-m-n}(x;d)=h_m(x;c)h_{M-m}(x;d)h_{n}(x;cq^m)
h_{N-n}(x;dq^{M-m}).$$
Performing the summation, using \eqref{tre}, gives
$$h_k(x;aq^{\alpha j})h_{M-k}(x;bq^{\beta(N-j)})
h_j(x;aq^{(1-\alpha)k})h_{N-j}(x;bq^{(1-\beta)(M-k)}).$$
For any  $\alpha,\beta\in\{0,1\}$, these factors combine to
$$h_{k+j}(x;a)h_{M+N-k-j}(x;b), $$
which completes the proof.
\end{proof}

\subsection{Combinatorial formulas}\label{cfss}

To get analogues of \eqref{mkc} and \eqref{kcomb}, we first consider
 all possible extensions of
 \eqref{hpf} to a general sum $h_{k_1+\dots+k_n}(x;a)$.
These are naturally labelled by permutations $\sigma$ of 
$\{1,\dots,n\}$:
$$h_{k_1+\dots+k_n}(x;a)=h_{k_{\sigma(1)}}(x;a)h_{k_{\sigma(2)}}
(x;aq^{k_{\sigma(1)}})\dotsm 
h_{k_{\sigma(n)}}(x;aq^{k_{\sigma(1)}+\dots+k_{\sigma(n-1)}}). $$
Replacing $\sigma$ by $\sigma^{-1}$, this may be written
$$h_{k_1+\dots+k_n}(x;a)=\prod_{i=1}^n
h_{k_i}(x;aq^{|k|_i^\sigma}),$$
where we introduced the notation
$$|k|_i^\sigma=\sum_{\{j;\,\sigma(j)<\sigma(i)\}} k_j $$
for a multi-index $k$. Note that
\begin{equation}\label{il}|k|_i^{\id}=k_1+k_2+\dots+k_{i-1}.\end{equation}

Thus, we have an extension of \eqref{rconv} labelled by two
 permutations $\sigma$, $\tau$:
\begin{multline}\label{rmc}
R_{k_1+\dots+k_n}^l(a,b,c,d;M_1+\dots+M_n)\\
=\sum_{m_1+\dots+m_n=l}\,\prod_{i=1}^n
R_{k_i}^{m_i}(aq^{|k|_i^\sigma},
bq^{|M-k|_i^\tau},cq^{|m|_i^{\id}},dq^{|M-m|_i^{\id}};M_i),
\end{multline}
where $0\leq k_i, m_i\leq M_i$. (We could replace both occurrences of $\id$ in 
\eqref{rmc} by a third permutation $\lambda$, but the resulting
identity is immediately reduced to \eqref{rmc} by permuting the $m_i$.)
In particular, when $M_1=\dots =M_n=1$, one has
\begin{multline}\label{f}
R_{k_1+\dots+k_n}^l(a,b,c,d;n)\\
=\sum_{\substack{m_1+\dots+m_n=l\\0\leq m_i\leq 1}}\,\prod_{i=1}^n
R_{k_i}^{m_i}(aq^{|k|_i^\sigma},
bq^{|1- k|_i^\tau},cq^{|m|_i^{\id}},dq^{|1- m|_i^{\id}};1).
\end{multline}

Note that on the right-hand side of \eqref{f}, only the elementary coefficients
$$R_k^m=R_k^m(a,b,c,d;1) $$
given by
$$\left(\begin{matrix}
R_0^0 & R_0^1\\ R_1^0 & R_1^1
\end{matrix}\right)=\left(\begin{matrix}
\displaystyle\frac{(1-bc)(1-b/c)}{(1-dc)(1-d/c)} &
\displaystyle\frac{(1-bd)(1-b/d)}{(1-cd)(1-c/d)} \\[3.5mm]
\displaystyle\frac{(1-ac)(1-a/c)}{(1-dc)(1-d/c)} &
\displaystyle\frac{(1-ad)(1-a/d)}{(1-cd)(1-c/d)}
\end{matrix}\right)$$
appear.
We shall see in Section \ref{ess} that the equation \eqref{f}
is closely related to the fusion of $R$-matrices
developped in \cite{d2,d3}. 
This explains the relation between our construction
and the statistical mechanics approach.

The combinatorics of the sum \eqref{f} deserves a separate
study, but we will make some further comments here. 
Note that, in \eqref{f}, a large number of right-hand sides give the
same left-hand side. If we only strive for a combinatorial
understanding of the coefficients $R_k^l$, it may be enough to choose
the right-hand side in a particularly simple fashion. For instance, we
may take
$\sigma=\tau=\id$, 
 and  choose $k_i$ as
\begin{equation}\label{k}
(k_1,\dots,k_n)=(\underbrace{1,\dots,1}_{k},\underbrace{0,\dots,0}_{n-k}).
\end{equation}

It seems natural to
 identify the summation-indices $m$ with lattice paths starting at
$(0,0)$ and going right at step $i$ if $m_i=1$ and up if $m_i=0$, thus
ending at $(l,n-l)$. Suppose that the $i$:th step in the path starts
at $(x,y)$. Then, by \eqref{il},
 $$|m|_i^{\id}=x,\qquad
|1- m|_i^{\id}=y. $$
Moreover, if the $k_i$ are chosen as in  \eqref{k}, then 
$$|k|_i^{\id}=\begin{cases}
i-1=x+y, & 1\leq i\leq k,\\
k, & k+1\leq i\leq n,
\end{cases}
$$
$$|1- k|_i^{\id}=\begin{cases}
0, & 1\leq i\leq k,\\
i-1-k=x+y-k, & k+1\leq i\leq n.
\end{cases}
$$
Thus, for instance, any one of the first $k$ steps in the path that
goes right  contributes a factor 
$$R_1^1(aq^{x+y},b,cq^x,dq^y;1)=\frac{(1-q^{x+2y}ad)(1-q^xa/d)}
{(1-q^{x+y}cd)(1-q^{x-y}c/d)}$$
to the sum. There are three other types of steps, giving rise to
similar factors. After replacing $n$ by $N$, this yields the following result.

\begin{corollary}\label{cc}
The coefficient $R_k^l(a,b,c,d;N;q)$ is given by the combinatorial
formula
\begin{multline*}
\sum_{\text{\emph{paths}}}\,
\prod_{\text{\emph{early right}}}\frac{(1-q^{x+2y}ad)(1-q^xa/d)}
{(1-q^{x+y}cd)(1-q^{x-y}c/d)}
\prod_{\text{\emph{early
up}}}\frac{(1-q^{2x+y}ac)(1-q^ya/c)}
{(1-q^{x+y}cd)(1-q^{y-x}d/c)}\\
\times\prod_{\text{\emph{late right}}}
\frac{(1-q^{x+2y-k}bd)(1-q^{x-k}b/d)}{(1-q^{x+y}cd)(1-q^{x-y}c/d)}
\prod_{\text{\emph{late
up}}}\frac{(1-q^{2x+y-k}bc)(1-q^{y-k}b/c)}{(1-q^{x+y}cd)(1-q^{y-x}d/c)}, 
\end{multline*}
where the sum is over all up-right lattice paths from $(0,0)$ to
$(l,N-l)$, the products are over steps in these paths,  the first $k$
steps being called ``early'' and the remaining $N-k$ steps being called 
``late''. In each factor, $(x,y)$ denotes the starting point of the
corresponding step.
\end{corollary}

There are many limit cases when Corollary \ref{cc} takes a
simpler form. It might be interesting to investigate  the limit cases
corresponding to various polynomials in the Askey Scheme. As an
example, let us consider the limit
$$L=\lim_{s\rightarrow 0}\lim_{c\rightarrow 0}(q^kd/bs)^{N-l}
R_k^l(as,bs,c,d;N;q). $$
It is easy to see from Theorem \ref{ret} that 
$$L=\qb{N}{l}(q^{k-l}a/b)^{k}\,
{}_{3}{\phi}_1\!\left[\begin{matrix} q^{-k},q^{-l},q^{k-N}a/b\\q^{-N}
\end{matrix};q,\frac{q^lb}{a}\right], $$
which, by \cite[Exercise 1.15]{gr} equals
\begin{equation}\label{qk}\qb{N}{l}(q^ka/b)^{k}\,
{}_{3}{\phi}_2\!\left[\begin{matrix} q^{-k},q^{-l},q^{-k}b/a\\q^{-N},0
\end{matrix};q,q\right].\end{equation}
(As an alternative, one may first use the
symmetries \eqref{symm} to write
$$R_k^l(a,b,c,d;N)=q^{-2\binom k2}a^{-2k}
q^{-2\binom{N-k}2}b^{-2(N-k)}R_k^l(q^{1-k}/a,q^{1+k-N}/b,c,d), $$
and then take the termwise limit in the accordingly transformed
version of Theorem \ref{ret}, thereby obtaining \eqref{qk} directly.)
The quantity \eqref{qk}
 may be identified with a $q$-Krawtchouk or dual $q$-Krawtchouk
polynomial \cite{ks}.
On the other hand, starting from Corollary \ref{cc} gives the
combinatorial expression
$$
L=\sum_{\text{{paths}}}\,
\prod_{\text{{early right}}} 1
\prod_{\text{{early
up}}}\frac{aq^{k+x}}{b}\prod_{\text{{late right}}}1
\prod_{\text{{late
up}}}q^x=\sum_{\text{{paths}}}\,\prod_{\text{{up}}}q^x
\prod_{\text{{early up}}}\frac{aq^k}{b}.
$$
Note that $\prod_{\text{{up}}}q^x=q^{\|\lambda\|}$, where
$\|\lambda\|$ is the number of boxes in  the Young
diagram to the upper left of the path.
Writing $t=aq^k/b$, we  conclude that
$$\qb{N}{l}t^{k}\,
{}_{3}{\phi}_2\!\left[\begin{matrix} q^{-k},q^{-l},1/t\\q^{-N},0
\end{matrix};q,q\right]
= \sum_{\text{{paths}}}q^{\|\lambda\|} t^{y(k)},
 $$
where $y(k)$ is the number of early ups, that is, 
 the $y$-coordinate of the end-point of the $k$:th step. 
This is a simple $q$-analogue of \eqref{kcomb}. Like \eqref{kcomb}, it
 is not very deep, but it gives an idea about what kind of information
 is contained in \eqref{f}. Note also that
 when $k=0$ or $t=1$ we recover the well-known fact
$$\qb{N}{l}=\sum_{\text{{paths}}}q^{\|\lambda\|}.$$

\section{Elliptic $6j$-symbols}\label{sell}

\subsection{Definition and elementary properties}
In this section we  discuss the extension of our approach to
 \emph{elliptic}
$6j$-symbols, or, more precisely, to their continuation in the parameters
studied in \cite{sz}. Roughly speaking, this corresponds to
 replacing everywhere ``$1-x$'' with the theta function
$$
\theta(x;p)=\prod_{j=0}^\infty(1-p^jx)(1-p^{j+1}/x),\qquad |p|<1.
$$
Since $\theta(x;0)=1-x$, the case
 $p=0$  will give back Wilson's functions discussed above.
The main difference is that in the elliptic case there is no
 Askey-type scheme of degenerate cases; all such  limits  require
$p=0$ to make sense.

We recall the notation \cite{gr2}
$$(a;q,p)_k=\prod_{j=0}^{k-1}\theta(aq^j;p), $$ 
$$\theta(x_1,\dots,x_n;p)=\theta(x_1;p)\dotsm\theta(x_n;p),$$
$$(a_1,\dots,a_n;q,p)_k=(a_1;q,p)_k\dotsm(a_n;q,p)_k. $$ 
Elliptic $6j$-symbols may be expressed in terms of the sum
\cite{gr2}
\begin{multline}\label{vd} {}_{12}V_{11}(a;b,c,d,e,f,g,q^{-n};q,p)\\
=\sum_{k=0}^n
\frac{\theta(aq^{2k})}{\theta(a)}\frac{(a,b,c,d,e,f,g,q^{-n};q,p)_k}
{(q,aq/b,aq/c,aq/d,aq/e,aq/f,aq/g,aq^{n+1};q,p)_k}\,q^k,
\end{multline}
subject to the balancing condition
$a^3q^{n+2}=bcdefg$. We mention that this function is
invariant under a natural action of $\mathrm{SL}(2,\mathbb Z)$
on $(q,p)$-space, cf.\ \cite{ft2,ths}.

Since 
\begin{equation}\label{ti}\theta(1/x;p)=-\theta(x;p)/x,\end{equation}
 the symbols $(a;q,p)_n$ satisfy
elementary identities similar to \eqref{eli}. Moreover,
\eqref{tadd} has the elliptic analogue (Riemann's addition formula)
\begin{equation}\label{thadd}\frac
vx\,\theta(xy,x/y,uv,u/v;p)=\theta(ux,u/x,vy,v/y;p)-\theta(uy,u/y,vx,v/x;p).
\end{equation} 

As an extension of \eqref{awm} we introduce the function
$$h_k(x;a)=h_k(x;a;q,p)=(a\xi,a\xi^{-1};q,p)_k,\qquad
\xi+\xi^{-1}=x.$$
For $a\neq 0$, this is an entire function of $x$.
 (If $a=0$, it does not make sense unless  $p=0$.)
We may then  introduce the
coefficients $R_k^l=R_k^l(a,b,c,d;N;q,p)$ by 
\begin{equation}\label{ee}
h_k(x;a)h_{N-k}(x;b)=\sum_{l=0}^N R_k^l(a,b,c,d;N;q,p)\,
 h_l(x;c)h_{N-l}(x;d).\end{equation}
Since the computation leading to Theorem \ref{ret} only used results
that have \emph{verbatim} elliptic extensions, it immediately carries
 over to the elliptic case.

\begin{theorem}\label{eet}
 For generic values of the parameters, the 
coefficients $R_k^l$ in \eqref{ee} exist uniquely and are given by
\begin{multline*}R_k^l(a,b,c,d;N;q,p)=
q^{l(l-N)}\frac{(q;q,p)_N}{(q;q,p)_l(q;q,p)_{N-l}}\\
\begin{split}&\times
\frac{(ac,a/c;q,p)_k(q^{N-l}bd,b/d;q,p)_l(b/c;q,p)_{N-k}(b/c;q,p)_{N-l}
(bc;q,p)_{N-k}}{(q^{l-N}c/d;q,p)_l(q^{-l}d/c;q,p)_{N-l}(cd;q,p)_N(b/c;q,p)_N
(bc;q,p)_l}\\
&\times{}_{12}V_{11}
(q^{-N}c/b;q^{-k},q^{-l},q^{k-N}a/b,q^{l-N}c/d,cd,q^{1-N}/ab,qc/b;q,p).
\end{split}\end{multline*}
\end{theorem}

\begin{remark}\label{tfr}
Like for Theorem \ref{ret},
the existence and uniqueness  falls out of the computation, but can
also be explained directly. Let $f$ be
any function of the form
\begin{equation}\label{mtf}
f(\xi)=\prod_{j=1}^{N}\theta(a_j\xi,a_j\xi^{-1};p),\end{equation}
and let $F$ be the function 
$$F(x)=f(e^{2\pi ix})=
\prod_{j=1}^{N}\theta(a_je^{2\pi ix},a_je^{-2\pi ix};p).$$
Then $F$ is an entire function satisfying
$$F(x+1)=F(x),\qquad F(x+\tau)=e^{-2\pi i N(2x+\tau)}F(x),\qquad F(-x)=F(x),$$
where $p=e^{2\pi i\tau}$. In classical terminology, $F$ is an even
theta function of order $2N$ and zero characteristics.  
It is known  that the space
$V_N$ of such functions has dimension $N+1$. We will denote by $W_N$
the space of corresponding functions $f$, that is, of
 holomorphic functions on $\mathbb
C\setminus\{0\}$ 
such that $f(\xi)=f(\xi^{-1})$ and
$f(p\xi)=(1/p\xi^2)^Nf(\xi)$. (In \cite{rai1}, these are called 
$BC_1$ theta functions of degree $N$.)
Now we observe that, for
$ab\neq 0$,  the functions
$$f_k(\xi)=h_k(x;a)h_{N-k}(x;b), \qquad x=\xi+\xi^{-1},$$
 are of the form \eqref{mtf}. With
essentially the same proof as for Lemma \ref{bl}, one may check that for 
$$p^ma/b\notin\{q^{1-N},q^{2-N},\dots,q^{N-1}\},\quad 
p^mab\notin\{1,q^{-1},\dots,q^{1-N}\},\qquad m\in\mathbb Z,$$
$(f_k)_{k=0}^N$ form a basis for $W_N$. We may then interpret 
Theorem \ref{eet} as giving the matrix for a change between two such
bases.
 \end{remark}

The following Corollary will be used below (in Section \ref{ess} and in
the proof of Proposition \ref{sal}).

\begin{corollary}\label{absc}
If $m$ and $n$ are non-negative integers, then
$$
h_k(x;a)h_{N-k}(x;b)\in\operatorname{span}_{k-m\leq l\leq k+n}\{
h_l(x;aq^m)h_{N-l}(x;bq^n)\}.$$
\end{corollary}

\begin{proof}
We need to consider the coefficient
$R_k^l(a,b,aq^m,bq^n;N;q,p)$. Isolating the factors containing the
quotients $b/d$ and $a/c$ in Theorem \ref{eet} gives
$$R_k^l(a,b,c,d;N;q,p)=\sum_{j=0}^{\min(k,l)}
\lambda_j (a/c;q,p)_{k-j}(b/d;q,p)_{l-j},$$
where $\lambda_j$ collects all other factors. If $a/c=q^{-m}$ this
vanishes unless $k-j\leq m$, and if $b/d=q^{-n}$
unless $l-j\leq n$. Thus, the range of summation is restricted to
$\max(k-m,l-n)\leq j\leq (k,l)$, which is empty unless  $k-m\leq l\leq k+n$.
This completes the proof.
\end{proof}

It is clear that the  coefficients $R_k^l$ enjoy similar
properties as were
obtained above in the special case
 $p=0$. This applies to the
biorthogonality relation \eqref{tbo}, the addition formula
\eqref{radd}, the convolution formulas in Corollary \ref{coco} and the
combinatorial formulas \eqref{f}. In particular, the identity in
Corollary \ref{coco} holds after replacing all factors $1-x$ with the
theta function $\theta(x;p)$.

\begin{remark}
Rains \cite{rai1} has obtained multivariable
 (Koornwinder--Macdonald-type)  extensions of the 
biorthogonal rational functions appearing in Theorem \ref{eet}.
The approach is  different from ours, although there
are similarities. Note, for instance, that the one-variable case
of the \emph{interpolation functions} in \cite[Definition 5]{rai1} are
essentially of the form $h_k(x;a)/h_k(x;b)$. 
 \end{remark}

\subsection{Comparison with statistical mechanics}\label{ess}

In this section we compare  the coefficients $R_k^l$ with  the
elliptic $6j$-symbols as defined in \cite{d3}. We shall see that
the latter correspond to certain discrete restrictions on
the parameters in $R_k^l$. 

We will 
follow the notation of \cite{d2}, where elliptic $6j$-symbols are denoted
$$W_{MN}(a,b,c,d|u).$$
 They depend on four external parameters $p$,
$\lambda$, $\xi$, $K$ and are defined for integers $a,b,c,d$ such that
\begin{gather}\label{wpc}
a-b,\ c-d\in\{-M,2-M,\dots,M\},\\ 
\label{wpc2} a-d,\ b-c\in\{-N,2-N,\dots,N\}.
\end{gather}
As was observed in \cite{ft2}, these symbols may be expressed
 in terms of the elliptic hypergeometric series
${}_{12}V_{11}$. Using Theorem \ref{eet}, we may then relate them to
the coefficients $R_k^l$. For instance, we have
\begin{multline}\label{wr}
W_{MN}(j+2l-N,i+2k-N,i,j|u)\\
=q^{(\xi+j+k+l-N)(l-k)+\frac 12
N(N-u)+\frac 12(N-2k)(i-j)
}\frac{(q^{u+M+1-N};q,p)_N}{(q;q,p)_N}\, R_k^l(a,b,c,d;N;q,p),
\end{multline}
where $q=e^{\pi i\lambda/K}$ and
\begin{multline}\label{abcd}
(a,b,c,d)\\
=(q^{\frac 12(u+\xi+i+1-N)},q^{\frac 12 (u-\xi-i+1-N)},
q^{\frac 12(u+\xi+j+M+1-N)},q^{\frac 12(u-\xi-j+M+1-N)}),\end{multline}
or, equivalently,
$$(q^M,q^{\xi+i},q^{\xi+j},q^u)=(cd/ab,a/b,c/d,abq^{N-1}) $$
($\xi$ is a parameter from \cite{d2} that has nothing to do with
\eqref{xi}). Note that \eqref{wpc2} corresponds to the 
condition $0\leq k,l\leq N$ on $R_k^l$, while \eqref{wpc} gives
 a further discrete restriction on the parameters.  

 In view of the large symmetry group of the terminating
${}_{12}V_{11}$, there are many different ways to identify elliptic
$6j$-symbols with the coefficients $R_k^l$. We have chosen the
representation \eqref{wr} 
since it explains the relation between fusion of $R$-matrices
 and the  combinatorial formulas of
Section \ref{cfss}. Namely, 
it is straight-forward to check that if we let $\sigma=\tau=\id$ in
\eqref{f}, specialize the parameters as in \eqref{abcd} and
replace $i$ by $n+1-i$ in the product, then \eqref{f} reduces
to \cite[Equation (2.1.21)]{d2}.

It is also interesting to consider the
 degeneration of the expansion problem
\eqref{ee} corresponding to the restriction \eqref{wpc}.
For this we introduce the parameter $m=(M+j-i)/2$.
The condition on $c-d$ in \eqref{wpc} means that $m$ is an integer with
$0\leq m\leq M$. 
The coefficients \eqref{wr} appear in the
expansion problem 
\begin{equation}\label{sjep}h_k(x;a)h_{N-k}(x;b)=\sum_{l=0}^N
R_k^l\,h_l(x;aq^m)h_{N-l}(x;bq^{M-m}).\end{equation}
By Corollary \ref{absc},  $R_k^l$ 
 vanishes unless $k-m\leq l\leq M+k-m$, which
corresponds exactly to the condition on $a-b$ in \eqref{wpc}. Then 
\eqref{sjep} reduces to 
$$h_k(x;a)h_{N-k}(x;b)=\sum_{l=\max(0,k-m)}^{\min(N,M+k-m)}
R_k^l\,h_l(x;aq^m)h_{N-l}(x;bq^{M-m}),\qquad 0\leq m\leq M, $$
which is thus the expansion problem solved by
 the elliptic $6j$-symbols of \cite{d3}.

\section{Sklyanin algebra and generalized eigenvalue problem}\label{ssa}

As was explained in Section \ref{rss}, our approach was motivated by 
previous work on  relations between the standard $\mathrm{SL}(2)$
quantum group    and  quantum $6j$-symbols. 
It is natural to ask what ``quantum group'' is behind the
more general case of elliptic $6j$-symbols. The answer turns out to be 
very satisfactory, namely, the  Sklyanin algebra \cite{sk1}.

\begin{figure}[htb]
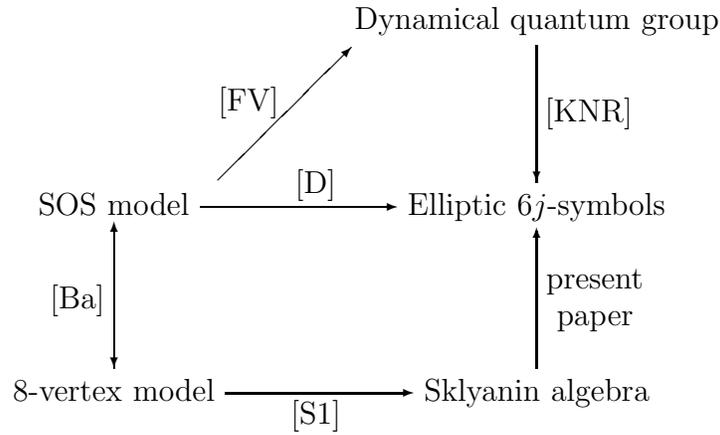
\label{qgc}
\begin{diagram}
& & \text{Dynamical quantum group}\\
&\ruTo^{\text{\cite{fv}}} & \dTo>{\text{\cite{knr}}} \\
\text{SOS model} & \rTo^{\quad\quad\text{\cite{d3}}}  &
\text{Elliptic $6j$-symbols} \\
\dCorresponds<{\text{\cite{b2}}}&
&  \uTo>{\begin{matrix}\text{present}\\\text{paper}\end{matrix}} \\
$8$\text{-vertex model} & \rTo_{\quad\quad
\text{\cite{sk1}}} & \text{Sklyanin algebra} 
\end{diagram}
\caption{Connections to quantum groups and solvable models}
\end{figure}

We recall that the Sklyanin algebra was  obtained from the $R$-matrix
of the eight-vertex model. Baxter found that this model is
related to a certain SOS (or face) model  by a
vertex-IRF transformation \cite{b2}. The original construction of
elliptic $6j$-symbols starts from the  $R$-matrix of the latter model.
Moreover, starting from Baxter's SOS model, Felder and Varchenko
constructed a dynamical quantum group \cite{fv}, which
was recently related to elliptic $6j$-symbols \cite{knr}.
We summarize these connections in Figure \ref{qgc}.

It would be interesting to find a direct link  between the 
approach of \cite{knr} and the discussion below. Presumably, this
 would involve  extending Stokman's paper \cite{st} to elliptic quantum
 groups. In particular, vertex-IRF transformations should play an
 important role.

To explain the connection with the Sklyanin algebra we introduce the
difference operators
\begin{multline*}
\Delta(a,b,c,d)f(\xi)\\
=\frac{\xi^{-2}\theta(a\xi,b\xi,c\xi,d\xi;p)
f(q^{\frac 12}\xi)- \xi^2\theta(a\xi^{-1},b\xi^{-1},c\xi^{-1},d\xi^{-1};p)
f(q^{-\frac 12}\xi)}{\xi\theta(\xi^{-2};p)}.\end{multline*}
Moreover, $N$ being fixed we write
$$\Delta(a,b,c)=\Delta(a,b,c,q^{-N}/abc).$$

The following observation was communicated to us by Eric Rains, see \cite{rai}.

\begin{proposition}[Rains, Sklyanin] \label{rsp}
The operators $\Delta(a,b,c)$
preserve the space $W_N$ defined in  \emph{Remark \ref{tfr}}.
 Moreover, they generate a
representation of the Sklyanin algebra on that space.
\end{proposition}

These representations were found by Sklyanin \cite[Theorem 4]{sk2},
except that he used the equivalent space denoted  $V_N$  in  Remark \ref{tfr}
and $\Theta_{00}^{2N+}$ in \cite{sk2}. Let $\Delta_i$, $i=0,1,2,3$
 be the operators representing Sklyanin's
generators $S_i$, pulled over from $V_N$ to $W_N$.  
Rains observed that every $\Delta_i$ is given
by an operator of the form $\Delta(a,b,c)$, with specific choices of
the parameters and that, conversely, every $\Delta(a,b,c)$ may be
expressed as a linear combination of the $\Delta_i$. 
One may view the resulting representation  as an
elliptic deformation of the group action \eqref{gr}.

Next we consider the action of the operators $\Delta$ on our
 basis vectors.

\begin{proposition}\label{sal}
With $x=\xi+\xi^{-1}$ one has
\begin{multline}\label{saa}
\Delta(a,b,c)h_k(x;q^{\frac 12}a)h_{N-k}(x;q^{\frac 12}b)\\
=\frac{q^{-N}}{abc}\,
\theta(q^kac,q^{N-k}bc,q^{N}ab;p)\,h_k(x;a)h_{N-k}(x;b)\end{multline}
and
\begin{equation}\label{sab}
\Delta(a,b,c)h_k(x;\lambda,\mu)\in\operatorname{span}_{k-1\leq j\leq
k+1}\{h_j(x;q^{\frac 12}\lambda,q^{\frac12}\mu)\}.\end{equation}
\end{proposition}

The identity \eqref{saa} is an analogue of the fact that, in the 
situation of Section \ref{sk}, any basis 
$((ax+b)^k(cx+d)^{N-k})_{k=0}^N$ is the eigenbasis of a Lie algebra
element. Similarly, \eqref{sab} is an analogue of the fact that any
other Lie algebra element acts tridiagonally on that basis.
The parameter shifts are unavoidable and related to the fact that
elliptic $6j$-symbols are biorthogonal rational functions rather than
orthogonal polynomials, see Remark \ref{ger} below.

\begin{remark}
The identities \eqref{saa} and \eqref{sab} are consistent in view of
$$
h_k(x;a)h_{N-k}(x;b)\in\operatorname{span}_{k-1\leq j\leq
k+1}\{h_j(x;aq)h_{N-j}(x;bq)\},$$
which is a special case of Corollary \ref{absc}. 
\end{remark}

Henceforth we suppress the deformation parameters $p$, $q$, thus
writing
$$\theta(x)=\theta(x;p),\qquad (a)_k=(a;q,p)_k.$$
When using notation such as $\theta(a\xi^{\pm})$, we will mean
$\theta(a\xi;p)\theta(a\xi^{-1};p)$. The following theta function
identity will be used in the proof of Proposition
\ref{sal}.

\begin{lemma}\label{pfl}
If $a_1\dotsm a_nb_1\dotsm b_{n+2}=1$, then
\begin{multline*}\xi^{-n-1}\prod_{j=1}^{n}\theta(a_j\xi)\prod_{j=1}^{n+2}
\theta(b_j\xi)-\xi^{n+1} \prod_{j=1}^n\theta(a_j\xi^{-1})
\prod_{j=1}^{n+2}\theta(b_j\xi^{-1})\\
=\frac{(-1)^n\xi\theta(\xi^{-2})}{a_1\dotsm
a_n}\sum_{k=1}^n\frac{\prod_{j=1}^{n+2}\theta(a_kb_j) \prod_{j=1,j\neq
k}^n\theta(a_j\xi^{\pm})}{\prod_{j=1,j\neq k}^n\theta(a_k/a_j)}. 
 \end{multline*} 
\end{lemma}

\begin{proof}
This is equivalent to the classical identity \cite[p.\
34]{tm}, see also \cite{r2},
$$\sum_{k=1}^n\frac{\prod_{j=1}^n\theta(a_k/b_j)}{\prod_{j=1,j\neq
k}^n\theta(a_k/a_j)}=0, \qquad a_1\dotsm a_n=b_1\dotsm b_n. $$
Namely, replace $n$ with $n+2$ and $b_j$ with $b_j^{-1}$ in that
identity, and put $a_{n+1}=\xi$, $a_{n+2}=\xi^{-1}$. Moving the last
two terms in the sum to the right gives
$$\sum_{k=1}^n\frac{\prod_{j=1}^{n+2}\theta(a_kb_j)}{\theta(a_k\xi^{\pm})\prod_{j=1,j\neq
k}^n\theta(a_k/a_j)}=-\frac{\prod_{j=1}^{n+2}\theta(\xi b_j)}{\theta(\xi^2)\prod_{j=1}^n\theta(\xi/a_j)}-\frac{\prod_{j=1}^{n+2}\theta(\xi^{-1}b_j)}{\theta(\xi^{-2})\prod_{j=1}^n\theta(\xi^{-1}/a_j)}.$$
After multiplying with $(-1)^n\xi\theta(\xi^{-2})\prod_{j=1}^n
a_j^{-1}\theta(a_j\xi^{\pm})$ and using \eqref{ti} 
repeatedly, one obtains the desired identity. 
\end{proof}

\begin{proof}[Proof of \emph{Proposition \ref{sal}}] We start with
\eqref{sab}. Writing out the left-hand side explicitly, collecting
common factors and using \eqref{ti} repeatedly gives
\begin{multline}\label{dah}
\Delta(a,b,c,d)\left((\lambda\xi^\pm)_k(\mu\xi^{\pm})_{N-k}\right)\\
\begin{split}&
=\frac{1}{\xi\theta(\xi^{-2})}\left\{\xi^{-2}\theta(a\xi,b\xi,c\xi,d\xi)\,
(q^{\frac 12}\lambda \xi, q^{-\frac12}\lambda\xi^{-1})_k( q^{\frac12}\mu\xi,q^{-\frac 12}\mu
\xi^{-1})_{N-k}\right.\\
&\left.\quad -\xi^2 \theta(a\xi^{-1},b\xi^{-1},c\xi^{-1},d\xi^{-1})\,
( q^{-\frac 12}\lambda\xi, q^{\frac12}\lambda\xi^{-1})_k(q^{-\frac12}\mu
\xi, q^{\frac 12}\mu\xi^{-1})_{N-k}
\right\}\\
&=\frac{q^{-1}\lambda\mu( q^{\frac 12}\lambda\xi^{\pm})_{k-1}
( q^{\frac 12}\mu\xi^{\pm})_{N-k-1}}{\xi\theta(\xi^{-2})}\\
&\quad\times\left\{\xi^{-4}\theta(a\xi,b\xi,c\xi,d\xi, q^{k-\frac
12}\lambda\xi,q^{\frac 12}\lambda^{-1}\xi,q^{N-k-\frac
12}\mu \xi,q^{\frac 12}\mu^{-1}\xi)\right.\\
&\left.\quad\quad -\xi^{4}\theta(a\xi^{-1},b\xi^{-1},c\xi^{-1},d\xi^{-1}, q^{k-\frac
12}\lambda\xi^{-1},q^{\frac 12}\lambda^{-1}\xi^{-1}, q^{N-k-\frac
12}\mu\xi^{-1},q^{\frac 12}\mu^{-1}\xi^{-1})\right\}.
\end{split}\end{multline}
Since $abcd=q^{-N}$,  we may  apply the case $n=3$ of Lemma \ref{pfl}
to the factor in brackets.
Choose  $(b_1,\dots,b_5)$ as $(a,b,c,d)$ together with any one
of the four numbers  
$$( q^{k-\frac 12}\lambda,q^{\frac 12}\lambda^{-1}, q^{N-k-\frac 12}\mu
,q^{\frac 12}\mu^{-1})$$
(let the audience pick it)
and choose $(a_1,a_2,a_3)$ as the remaining  three of those
 numbers. As a function of $\xi$, our expression then takes the form
\begin{multline*}( q^{\frac 12}\lambda\xi^{\pm})_{k-1}
(q^{\frac
12}\mu \xi^{\pm})_{N-k-1}\\
\times\left\{C_1\theta(a_2\xi^{\pm},a_3\xi^{\pm})+
C_2\theta(a_1\xi^{\pm},a_3\xi^{\pm})+C_3\theta(a_1\xi^{\pm},
a_2\xi^{\pm})\right\}.\end{multline*}
Depending on the choice of $a_i$, each  term is proportional
to one of the six functions
$$(q^{-\frac 12}\lambda \xi^{\pm})_{k}(q^{-\frac
12}\mu\xi^{\pm})_{N-k}, \quad (q^{-\frac 12}\lambda \xi^{\pm})_{k+1}(q^{\frac
12}\mu\xi^{\pm})_{N-k-1}, \quad (q^{-\frac 12}\lambda \xi^{\pm})_{k}(q^{\frac
12}\mu\xi^{\pm})_{N-k},$$
$$(q^{\frac 12}\lambda \xi^{\pm})_{k}(q^{-\frac
12}\mu\xi^{\pm})_{N-k}, \quad (q^{\frac 12}\lambda \xi^{\pm})_{k-1}(q^{-\frac
12}\mu\xi^{\pm})_{N-k+1}, \quad (q^{\frac 12}\lambda \xi^{\pm})_{k}(q^{\frac
12}\mu\xi^{\pm})_{N-k}.$$
By Corollary \ref{absc},  these  all belong to 
$$\operatorname{span}_{k-1\leq j\leq k+1}
\{(q^{\frac 12}\lambda\xi^{\pm})_j(q^{\frac 12}\mu\xi^\pm)_{N-j}\}. $$
This completes the proof of \eqref{sab}.

If we put $\lambda =q^{\frac 12}a$, $\mu=q^{\frac 12}b$ in
\eqref{dah}, the  factor $\theta(a\xi^{\pm},b\xi^{\pm})$ can
be pulled out from the bracket, giving
\begin{multline*}
\Delta(a,b,c,d)\left((q^{\frac  12}a\xi^\pm)_k(q^{\frac 12}b
\xi^{\pm})_{N-k}\right)
=\frac{(a\xi^{\pm})_{k}
(b\xi^{\pm})_{N-k}}{\xi\theta(\xi^{-2})}\\
\times
\left\{\xi^{-2}\theta(c\xi,d\xi, q^{k}a\xi, q^{N-k}b\xi)
 -\xi^{2}\theta(c\xi^{-1},d\xi^{-1}, q^{k}a\xi^{-1},
q^{N-k}b \xi^{-1})\right\}. 
\end{multline*}
The case $n=1$ of Lemma \ref{pfl}, which is
equivalent to \eqref{thadd}, now gives \eqref{saa}.
\end{proof}

\begin{remark}\label{ger}
Proposition \ref{sal}  connects our work with the generalized
eigenvalue problem (GEVP), which is central to the approach of
Spiridonov and Zheda\-nov \cite{sz,sz2}. Recall that, roughly speaking,
 the theory of orthogonal polynomials is equivalent to spectral
theory of Jacobi operators, that is, to the eigenvalue problem
$$Ye_k=\lambda_k e_k $$
for a (possibly infinite) tridiagonal matrix $Y$. The theory of
biorthogonal rational functions similarly corresponds to the GEVP
\begin{equation}\label{gevp}Y_1e_k=\lambda_k Y_2e_k\end{equation}
for two tridiagonal matrices $Y_1$, $Y_2$.
Note that \eqref{saa} means that  
$$e_k=h_k(x;a)h_{N-k}(x;b)$$
solves the two-parameter family of GEVP:s
$$\Delta_1e_k=\lambda_k\Delta_2e_k, $$
where $\Delta_1=\Delta(q^{-\frac 12}a,q^{-\frac 12}b,c)$,
$\Delta_2=\Delta(q^{-\frac 12}a,q^{-\frac 12}b,d)$, with 
 $c$ and $d$  arbitrary. Moreover, if we let $\Delta_3$ be any operator
of the form $\Delta(e,f,g)$ and we put $Y_1=\Delta_3\Delta_1$,
$Y_2=\Delta_3\Delta_2$, we have that $e_k$ solves \eqref{gevp} with
$Y_i$  tridiagonal in the basis $(e_k)_{k=0}^N$.
 Thus, we may view 
elliptic $6j$-symbols as the change of base matrix between the
solutions of two different GEVP:s, where the involved tridiagonal operators
 are appropriate elements of the
Sklyanin algebra, acting in a finite-dimensional representation.
\end{remark}

\end{document}